\documentclass[letterpaper, 11pt,  reqno]{amsart}

\usepackage{amsmath,amssymb,amscd,amsthm,amsxtra, esint}

\usepackage[implicit=true]{hyperref}

\allowdisplaybreaks[2]

\usepackage{color}

\definecolor{gr}{rgb}   {0.,   0.69,   0.23 }
\definecolor{bl}{rgb}   {0.,   0.5,   1. }
\definecolor{mg}{rgb}   {0.85,  0.,    0.85}
\definecolor{yl}{rgb}   {0.8,  0.7,   0.}
\definecolor{or}{rgb}  {0.7,0.2,0.2}

\newtheorem{theorem}{Theorem} [section]

\newtheorem{lemma}[theorem]{Lemma}
\newtheorem{proposition}[theorem]{Proposition}
\newtheorem{remark}[theorem]{Remark}

\newtheorem{corollary}[theorem]{Corollary}

\newtheorem*{ackno}{Acknowledgment}

\DeclareMathOperator*{\intt}{\int}

\DeclareMathOperator*{\supp}{supp}
\DeclareMathOperator{\med}{med}

\newcommand{\noi}{\noindent}
\newcommand{\Z}{\mathbb{Z}}
\newcommand{\R}{\mathbb{R}}

\newcommand{\T}{\mathbb{T}}

\newcommand{\sech}{\textup{sech}}

\let\Re=\undefined\DeclareMathOperator*{\Re}{Re}

\let\P= \undefined
\newcommand{\P}{\mathbf{P}}

\newcommand{\N}{\mathbb{N}}

\newcommand{\FL}{\mathcal{F}L} 

\renewcommand{\S}{\mathcal{S}}

\newcommand{\F}{\mathcal{F}}

\newcommand{\dl}{\delta}

\newcommand{\eps}{\varepsilon}

\newcommand{\ld}{\lambda}

\newcommand{\s}{\sigma}

\newcommand{\ft}{\widehat}

\newcommand{\wt}{\widetilde}
\newcommand{\cj}{\overline}
\newcommand{\dx}{\partial_x}
\newcommand{\dt}{\partial_t}

\newcommand{\ta}{\theta}

\renewcommand{\l}{\ell}

\newcommand{\les}{\lesssim}
\newcommand{\ges}{\gtrsim}

\newcommand{\jb}[1]
{\langle #1 \rangle}

\newcommand{\ind}{\mathbf 1}

\numberwithin{equation}{section}
\numberwithin{theorem}{section}

\begin{document}

\baselineskip = 14pt

\title
[GWP of  mKdV in modulation spaces]
{On global well-posedness of the modified KdV equation in modulation spaces}

\author[T.~Oh and  Y.~Wang]
{Tadahiro Oh and Yuzhao Wang}

\address{
Tadahiro Oh, School of Mathematics\\
The University of Edinburgh\\
and The Maxwell Institute for the Mathematical Sciences\\
James Clerk Maxwell Building\\
The King's Buildings\\
Peter Guthrie Tait Road\\
Edinburgh\\ 
EH9 3FD\\
 United Kingdom}

\email{hiro.oh@ed.ac.uk}

\address{
Yuzhao Wang\\
School of Mathematics\\
Watson Building\\
University of Birmingham\\
Edgbaston\\
Birmingham\\
B15 2TT\\ United Kingdom}

\email{y.wang.14@bham.ac.uk}



\begin{abstract}
We study well-posedness of the complex-valued modified KdV equation (mKdV)
on the real line.
In particular, we prove  local well-posedness of mKdV in modulation spaces $M^{2,p}_{s}(\R)$ for $s \ge \frac14$
and $2\leq p < \infty$.
For $s < \frac 14$, we show that the solution map for mKdV is not locally uniformly continuous
in $M^{2,p}_{s}(\R)$.
By combining this local well-posedness with our previous work (2018) on 
an a priori global-in-time bound for mKdV in modulation spaces, 
we also establish  global well-posedness of mKdV 
in $M^{2,p}_{s}(\R)$ for $s \ge \frac14$
and $2\leq p < \infty$.

\end{abstract}


\maketitle
%



\section{Introduction}

\subsection{Modified KdV equation}
We consider  the Cauchy problem for the  complex-valued modified KdV equation 
 on the real line:
\begin{align}
\begin{cases}
 \dt u + \dx^3 u \pm 6 |u|^2 \dx u =0 \\
u|_{t= 0} = u_0, 
\end{cases}
\qquad (x, t) \in \R\times \R.
\label{mKdV}
\end{align}

\noi
The equation \eqref{mKdV} is known to be completely integrable 
and is closely related to the cubic nonlinear Schr\"odinger equation (NLS):
\begin{align}
i \dt u -  \dx^2 u \mp 2 |u|^2 u = 0.
\label{NLS}
\end{align}

\noi
See  \cite{Hirota, SS, KVZ, OW}.
When the initial data $u_0$ is real-valued, 
the corresponding solution $u$ to \eqref{mKdV} remains real-valued, 
thus solving
the following real-valued mKdV:
\begin{align}
 \dt u + \dx^3 u \pm 6 u^2 \dx u = 0.
\label{mKdV2}
\end{align}

\noi
The mKdV enjoys the following scaling symmetry: 
\begin{align}
u(x, t) \longmapsto
u_\ld(x, t) = \ld^{-1} u (\ld^{-1}x, \ld^{-3}t), 
\label{scaling}
\end{align}

\noi
which induces the scaling-critical Sobolev regularity $s_\text{crit} = - \frac 12$
in the sense that 
homogeneous $\dot{H}^{-\frac 12}$-norm is invariant
under the scaling symmetry \eqref{scaling}.

The Cauchy problem \eqref{mKdV} has been studied extensively.
In \cite{Kato}, Kato studied \eqref{mKdV}
from a viewpoint of quasilinear hyperbolic equations
(in particular, not making use of dispersion)
and proved its local well-posedness 
in $H^s(\R)$, $s > \frac32$.
In \cite{KPV89, KPV93}, 
Kenig-Ponce-Vega 
exploited the dispersive nature of the equation
and proved local well-posedness of \eqref{mKdV} in $H^s(\R)$, $s \ge \frac14$.
In \cite{Tao}, 
Tao  gave an alternative proof of the local well-posedness in $H^{\frac14}(\R)$ 
by using the Fourier restriction norm method.
We also mention recent papers~\cite{MPV, KOY} on 
unconditional uniqueness of solutions to 
 \eqref{mKdV}
in $H^s(\R)$, $s > \frac 14$. 
Let us now turn our attention to global well-posedness
of \eqref{mKdV}.
In the real-valued setting, 
Colliander-Keel-Staffilani-Takaoka-Tao~\cite{CKSTT} applied the $I$-method
and proved global well-posedness of \eqref{mKdV2}
in $H^s(\R)$ for $s > \frac 14$.
See  Kishimoto~\cite{Kishimoto0} for the endpoint global well-posedness
 in $H^{\frac14}(\R)$.
In a recent paper, 
Killip-Vi\c{s}an-Zhang \cite{KVZ} exploited the 
completely integrable structure of the equation and 
proved a global-in-time a priori bound on the $H^s$-norm of solutions to 
the complex-valued mKdV \eqref{mKdV}
for $-\frac 12 < s < 0$.
While it is not written in an explicitly manner,\footnote{See also Appendix B of \cite{OW}
for details of a global-in-time a priori bound 
on the $H^s$-norm of solutions 
to the complex-valued mKdV \eqref{mKdV}
for $0 < s < \frac 12$.} 
their result is readily extendable
to $-\frac 12 < s < 1$
and thus proves global well-posedness of 
the complex-valued mKdV \eqref{mKdV} in $H^{\frac 14}(\R)$.

On the other hand, 
it is known that the solution map to \eqref{mKdV} is not locally uniformly continuous
in $H^s(\R)$ for  $s < \frac 14$; see \cite{KPV01, CCT1}.
This in particular 
implies that one can {\it not} use a contraction argument 
to construct solutions to \eqref{mKdV} in this regime.
One possible approach  to study rough solutions outside $H^\frac 14 (\R)$ is
to use a more robust energy method.
In \cite{CHT}, Christ-Holmer-Tataru 
employed an energy method in the form of the short-time Fourier restriction norm
method and proved global existence of solutions
to the real-valued mKdV~\eqref{mKdV2} in $H^s(\R)$ for $-\frac 18 < s < \frac 14$.
Uniqueness of these solutions, however, is unknown at this point.

Another approach is to study the Cauchy problem \eqref{mKdV}
in some other scales of function spaces than the Sobolev spaces $H^s(\R)$.
In \cite{Grun}, 
Gr\"unrock studied \eqref{mKdV}
in the Fourier-Lebesgue spaces $\FL^{s, p}(\R)$ defined by the norm:
\begin{align*}
\| f \|_{\FL^{s, p}} = \| \jb{\xi}^s \ft f (\xi)\|_{L^p},
\end{align*}

\noi
where $\jb{\,\cdot\,} = (1+|\cdot|^2)^\frac{1}{2}$,
and proved its local well-posedness
in $\FL^{s, p}(\R)$
for $s \geq \frac{1}{2p}$ and $2 \leq p < 4$.
In  \cite{GV}, Gr\"unrock-Vega extended this result 
to $2 \leq p < \infty$ with the same range of $s\geq \frac{1}{2p}$.
Note that 
the space $\FL^{0, \infty}(\R)$ of pseudo-measures
is critical
in terms of the scaling symmetry \eqref{scaling}, 
i.e.~the $\FL^{0, \infty}$-norm remains invariant under \eqref{scaling}.
Hence, by taking $p \to \infty$, 
we see that 
the local well-posedness result in \cite{GV} is almost critical.
There are two remarks in order;
(i) the range of $s \geq \frac{1}{2p}$
in \cite{Grun, GV} is sharp in the sense that 
the solution map to \eqref{mKdV} is not locally uniformly continuous
in $\FL^{s, p}(\R)$ for  $s < \frac 1{2p}$ and $2 \leq p < \infty$.
See Section 5 in \cite{Grun}.
(ii) there seems to be no known global well-posedness of \eqref{mKdV}
in the context of Fourier-Lebesgue spaces, 
extending local solutions constructed in \cite{Grun, GV} globally in time.

\subsection{Main results}
Our main goal in this paper is to study the Cauchy problem \eqref{mKdV}
in modulation spaces $M^{2, p}_s(\R)$.
We first 
 recall  the definition of modulation spaces $M^{r, p}_s(\R)$ introduced in \cite{FG1,FG2}.
Let $\psi \in \S(\R)$ such that
\begin{align*}
\supp \psi \subset [-1, 1]
\qquad \text{and} \qquad \sum_{k \in \Z} \psi(\xi -k) \equiv 1.
\end{align*}

\noi
Then, given $s\in\R$, $1\leq r, p \leq \infty$, the modulation space $M_s^{r, p}(\R)$ is defined
as the collection of all tempered distributions
$f\in\S'(\R)$ such that
$\|f\|_{M_s^{r, p}}<\infty$, where
the $M_s^{r, p}$-norm is defined by
\begin{equation*} 
\|f\|_{M_s^{r, p}} 
= \big\| \jb{n}^s\|\psi_n(D) f \|_{L_x^r(\R)} \big\|_{\l^p_n(\Z)}.
\end{equation*}

\noi
Here, 
$\psi_n(D)$ is the Fourier multiplier operator
with the multiplier 
\begin{align*}
\psi_n(\xi) := \psi(\xi - n).
\end{align*}

\noi
In the following, we only consider $r = 2$.
In the case of $r=2$, we have the following embedding 
\begin{align}
M^{2, p}_s(\R) \supset \FL^{s,p}(\R)
\label{mod3}
\end{align}

\noi
for $p\geq 2$.
The embedding \eqref{mod3} is immediate from
\[ \| f\|_{\FL^{s, p}} \sim 
\big\| \jb{n}^s\| \psi_n (\xi) \ft f (\xi)\|_{L_\xi^p(\R)} \big\|_{\l^p_n(\Z)}\]

\noi
and the support condition on $\psi$.

In \cite{OW}, we extended the work \cite{KVZ} by Killip-Vi\c{s}an-Zhang
on the global-in-time a priori bound
for solutions to \eqref{mKdV} to the modulation space setting
and proved the following result.

\begin{proposition}\label{PROP:global}
Let $2 \leq p < \infty$ and $0 \leq s < 1 - \frac 1p$.\footnote{The upper bound $1- \frac 1p$ is not essential
and we expect that this restriction  can be relaxed by a consideration similar to that in Section 3 of \cite{KVZ}.} 
Then, 
there exists $C= C(p) >0$ such that 
\begin{align}
\|u(t)\|_{M^{2, p}_s}
\leq C \big(\|u(0)\|_{M^{2, p}_s}\big)
\label{bd1}
\end{align}

\noi
for any Schwartz class  solution $u$ to the complex-valued  mKdV \eqref{mKdV}
and  any $t \in \R$.

\end{proposition}

In \cite{OW}, 
we also established  the same global-in-time a priori bound
for solutions to the cubic NLS \eqref{NLS}.
Combining this with the local well-posedness
of \eqref{NLS} in 
$M^{2, p}_s(\R)$ for $s \geq 0$ and $2 \leq p < \infty$
by S.~Guo \cite{Guo}, 
we proved global well-posedness of the cubic NLS \eqref{NLS}
in almost critical modulation spaces\footnote{The modulation spaces
 are based on the unit cube decomposition of the frequency space
 and thus there is no scaling for the modulation spaces.
 We, however, say that $M^{2, \infty}_0(\R)$ is a critical space
 in view of the embedding~\eqref{mod3}
 with $s  = 0$ and $p = \infty$.}
$M^{2, p}_s(\R)$ for $s \geq 0$

On the other hand, 
there is no known local well-posedness for the modified
KdV equation~\eqref{mKdV} in the modulation space $M^{2, p}_s(\R)$,
which motivated us to prove the following local well-posedness result.

\begin{theorem}
\label{THM:1}
Let $s \geq \frac 14$ and $2\le p<\infty$.
Then,
the complex-valued mKdV \eqref{mKdV}
is locally well-posed in $M^{2, p}_s(\R)$.
\end{theorem}

In \cite{Guo}, 
S.~Guo proved local well-posedness
of  the cubic NLS \eqref{NLS} in the modulation spaces $M_s^{2,p}(\R)$ 
for $s \geq 0$ and $2 \leq p < \infty$.
The proof was based on the Fourier restriction norm method
adapted to the modulation spaces, 
where  an endpoint version of two-dimensional 
Fourier restriction estimate played a crucial role.
See also \cite{GRW} for a work on  the derivative NLS 
which employs a similar strategy.
In proving Theorem \ref{THM:1}, 
we also use the Fourier restriction norm method adapted to the  modulation space setting.
See \eqref{Xsb} below.
We, however, provide a different approach than~\cite{Guo, GRW}.
Our argument is  based on bilinear estimates;  see  Lemmas \ref{LEM:bi-1} and 
Corollary~\ref{COR:bi}.
It is worthwhile to mention that our approach works equally well for the cubic NLS and the derivative NLS,
providing an  alternative approach to the  results in \cite{Guo,GRW}.

As a corollary to Proposition \ref{PROP:global}
and Theorem \ref{THM:1}, 
we obtain the following global well-posedness.

\begin{theorem}
\label{THM:2}
Let $s \geq \frac 14$ and $2\le p<\infty$.
Then,  there exists a function $C: \R_+ \times \R_+ \to \R_+$, which is increasing in each argument, 
such that 
\begin{align}
\sup_{t \in [-T, T]} \| u(t) \|_{M^{2, p}_s} \leq
 C\big(\|u_0\|_{M^{2, p}_s}, T\big)
\label{non5}
\end{align}

\noi
for any $T > 0$
and any Schwartz solution $u$ to \eqref{mKdV} with $u|_{t = 0} = u_0$.
In particular, this implies that the complex-valued mKdV \eqref{mKdV}
is globally well-posed in $M^{2, p}_s(\R)$.

\end{theorem}

For $\frac 14 \leq s < 1- \frac {1}{p}$,
Proposition \ref{PROP:global}
allows us to choose 
the right-hand side of \eqref{non5} to be independent of $T>0$.
For $s \geq 1- \frac{1}{p}$, 
we combine a persistence-of regularity argument
with the global-in-time bound on the $M^{2, p}_\frac{1}{4}$-norms
of solutions.  See Subsection \ref{SUBSEC:GWP}.

\begin{remark}\rm
One can easily adapt the proof of Theorem \ref{THM:1}
to extend the local well-posedness of \eqref{mKdV}
to  $1\le p<2$ (and   $s \geq \frac 14$).
Similarly, by establishing persistence of regularity
as in \cite{OW2}, 
we can also prove global well-posedness of \eqref{mKdV}
in $M^{2, p}_s(\R)$ 
for $s \geq \frac 14$ and $1\leq p < 2$.
See Remark \ref{REM:pers}.

\end{remark}

\smallskip

On the one hand, 
$\dot \FL^{\frac 14 , \infty}(\R)$
scales like $\dot H^{-\frac 14}(\R)$
and thus we may say that 
$ M^{2, \infty}_{\frac 14}(\R) $ ``scales like'' 
$\dot H^{-\frac 14}(\R)$ in view of the embedding \eqref{mod3}.
On the other hand, 
the $M^{2, p}_{s}(\R)$-norm is weaker
than the $\FL^{s, p}$-norm for $p > 2$
and the solution map to the mKdV \eqref{mKdV}
 fails to be locally uniformly continuous
in $M^{2, p}_{s}(\R)$ as soon as $s < \frac 14$.

\begin{proposition}
\label{PROP:ill}
Suppose that $(s, p)$ satisfies one of the following conditions:
\textup{(i)}  $2 \leq p  \leq  \infty$ and $0 \leq  s < \frac 14$
or \textup{(ii)} $2 \leq p  <  \infty$ and $  - \frac 1p < s < 0$.
%
Then,   the data-to-solution map for~\eqref{mKdV} 
in the focusing case \textup{(}with the $+$ sign in \eqref{mKdV}\textup{)\, :}
$u_0 \in M^{2, p}_s\mapsto u \in C([-T, T]; M^{2, p}_s(\R))$
is not locally uniformly continuous
for any $T > 0$.
\end{proposition}

Proposition \ref{PROP:ill} shows   a sharp contrast with the Fourier-Lebesgue case,
where local well-posedness was proved via a contraction argument even for some $s < \frac 14$.

In \cite{KPV01}, 
Kenig-Ponce-Vega proved the failure of local uniform continuity 
of the solution map for the complex-valued focusing mKdV~\eqref{mKdV} in $H^s(\R)$, $- \frac 12 < s < \frac 14$,
by building counterexamples from explicit soliton solutions.
See~\eqref{un} below.
By making use of breather solutions to the real-valued focusing mKdV \eqref{mKdV}, 
they also extended this  result for the real-valued case.
In~\cite{CCT1}, 
Christ-Colliander-Tao \cite{CCT1} extended this failure
of local uniform continuity below $H^\frac{1}{4}(\R)$ 
(for $-\frac 14 < s < \frac 14$) to the defocusing case
by approximating the mKdV dynamics by the cubic NLS dynamics
(which was in turn approximated by a dispersionless equation).
These (approximate) solutions in \cite{KPV01, CCT1}
depend on a parameter $N$ tending to $\infty$
and,   as $N \to  \infty$, 
they start to concentrate at a single point on the frequency side
(for $s > 0$).
Namely, they are essentially supported on a single unit cube
for $N \gg 1$.
In this regime, their $M^{2, p}_s$-norms
 basically reduce to  the $H^s$-norms, 
 giving rise to the threshold regularity $s = \frac 14$
 even in the modulation space setting.
 The main difficulty is that calculation required
 for the modulation space setting
is  much more involved than that for the Sobolev space setting.
Therefore, we only demonstrate the proof for the focusing cases in Section \ref{SEC:ill}.
We expect the same result hold for the defocusing case.
For the conciseness of the paper, however, 
we choose not to discuss details for the defocusing case.

\begin{remark}\rm

In a recent 
preprint \cite{CG}, the authors
independently proved  local well-posedness
of~\eqref{mKdV}
analogous to Theorem \ref{THM:1} 
for  $s \geq \frac 14$ and $2\le p\leq \infty$.
While the result in \cite{CG} only refers to local well-posedness,
it contains the $p = \infty$ case.
In view of the embedding 
$M^{2, \infty}_s  (\R) \subset M^{2,p }_\frac{1}{4}(\R)$
for 
\begin{align}
(s - \tfrac 14) p >1,
\label{non6} 
\end{align}
a combination of  the a priori bound \eqref{non5}
in Theorem \ref{THM:2}
(with $s > \frac 14$ and $p < \infty$ satisfying~\eqref{non6})
and a persistence-of-regularity argument
as in Subsection \ref{SUBSEC:GWP}
seems to yield global well-posedness of \eqref{mKdV}
for $s > \frac 14$ and $p = \infty$.
On the other hand, the global well-posedness issue
at the endpoint case: $s = \frac 14$ and $p = \infty$
remains open.

\end{remark}

\section{Preliminaries} 

Given dyadic $N\geq 1$, 
we denote by  $\P_N$  the Littlewood-Paley projector
onto the (spatial) frequencies $\{|\xi|\sim N\}$.
We use the following convention;
any summation over capitalized variables such as 
$N_1$, $N_2$, $\dots$, are presumed to be over dyadic numbers
of the form $2^k$, $k \in  \N \cup \{0\}$.

For $n \in \Z$, 
let 
\begin{align}
\ft{\Pi_n f}(\xi) = \psi_n(\xi) \ft f(\xi) = 
 \psi_n(\xi - n) \ft f(\xi).
\label{Pi1}
\end{align}

\noi
By Bernstein's inequality, we have
\begin{align}
\begin{split}
\| \P_N f\|_{L^p_x}
& \les N^{\frac{1}{q} - \frac{1}{p}} \| f \|_{L^q_x}, \\
\| \Pi_n f\|_{L^p_x}
& \les  \| f \|_{L^q_x} 
\end{split}
\label{bern}
\end{align}

\noi
for any $1 \leq q \leq p \leq  \infty$.

In the seminal work \cite{Bo93}, 
Bourgain introduced
the $X^{s, b}$-space defined by the norm:
\begin{align*}
\|u\|_{X^{s,b}} 
 : =  \| \jb{\xi}^s  \jb{\tau-\xi^3}^{b} \ft u (\xi,\tau) \|_{L_{\tau,\xi}^2}.
\end{align*}

\noi
In this paper, we use the following version of the $X^{s, b}$-space adapted
to the modulation spaces $M^{2, p}_s(\R)$:
\begin{align}
\|u\|_{X^{s,b}_{p}} 
& : = \bigg( \sum_{n\in \Z}  \jb{n}^{sp} \|  \jb{\tau-\xi^3}^{b} \ft u (\xi,\tau) \|_{L_{\tau,\xi}^2(\R\times [n,n+1])}^p  \bigg)^{\frac1p} \notag \\
& \sim \big\| \|\Pi_n  u \|_{X^{s,b}} \big\|_{\l^p_n }.
\label{Xsb}
\end{align}

\noi
When $p = 2$, the space $X^{s, b}_p$ reduces
to  the usual $X^{s, b}$-space.
When $b > \frac 12$, the following embedding holds:
\begin{align}
X^{s, b}_p \subset C(\R; M^{2, p}_s(\R)). 
\label{embed1}
\end{align}

Let $p \geq q \geq 1$.
Since $\l^q_n (\Z) \subset \l^p_n(\Z)$, we have 
\begin{align}
\|u\|_{X^{s,b}_{p}} \le \|u\|_{X^{s,b}_{q}}.
\label{embed1a}
\end{align}

\noi
On the other hand, 
from H\"older's inequality, we have 
\begin{align}
\|\P_N u\|_{X^{s,b}_{q}} \les N^{\frac1q - \frac1p}  \|\P_N u\|_{X^{s,b}_{p}}.
\label{embed2}
\end{align}


\noi

Given a time interval $I \subset \R$, 
 we also define the local-in-time version $X^{s,b}_{p}(I)$ of 
the  $X^{s, b}_{p}$-space
as the collection of functions $u$ such that 
\begin{equation*}
\|u\|_{X^{s,b}_{p}(I)}:=\inf \big\{\|v\|_{X^{s,b}_{p}} 
\, : \,  v|_{I}=u \big\}
\end{equation*}

\noi
is finite.

The following linear estimates follow
from  the characterization \eqref{Xsb}
and 
the corresponding linear estimates
for the standard $X^{s, b}$-spaces.
See \cite{GTV} for the proof.

\begin{lemma}\label{LEM:lin} 
\textup{(i) (Homogeneous linear estimate).}
Given $1\leq p< \infty$  and $s,b\in \R$, we have 
\begin{equation*}
\|e^{-t\dx^3} f\|_{X^{s,b}_{p}([0, T])} \les \|f\|_{M^{2, p}_s}
\end{equation*}

\noi
for any $0 < T \leq 1$.

\medskip

\noi
\textup{(ii) (Nonhomogeneous linear estimate).}
Let $s\in \R$, $1\leq p <  \infty$,  and $-\frac{1}{2}<b'\leq 0\leq b \leq 1+b'.$ 
Then, we have 
\begin{equation*}
\bigg\| \int_0^te^{-(t-t')\dx^3} F(t')dt'\bigg\|_{X^{s,b}_{p}([0, T])}\les T^{1+b'-b}\|F\|_{X^{s,b'}_{p}([0, T])}
\end{equation*}

\noi
for any $0< T \leq 1$.

\end{lemma}

In the following, we list 
various estimates in proving the crucial trilinear estimate
(Proposition \ref{PROP:trilinear}).
The following inequality will be convenient
in dealing with the resonant case in Section \ref{SEC:tri}.
From  H\"older's and Young's inequalities, we have
\begin{align}
\label{hlp2}
\sum_{\substack{m, n \in \Z\\ m\neq n}} \frac{a_m b_n}{|m-n|\jb{n}^\eps } \leq C_\eps   \|a_n\|_{\l^p(\Z)} \|b_n\|_{\l^{p'}(\Z)}
\end{align}

\noi
for any $\eps > 0$, 
where $p'$ denotes the H\"older conjugate of $p$.

Next, we recall  a bilinear estimate
from \cite{Grun}.  
Given $\ta > 0$, let $I^\ta=(-\dx^2)^\frac{\ta}{2}$
denote the Riesz potential of order $\theta$. 
We also define $I^\ta_-$
by 
\[
\F_x (I^\theta_-(f,g))(\xi) : = \int_{\xi = \xi_1+\xi_2} |\xi_1-\xi_2|^\theta \ft f (\xi_1) \ft g(\xi_2) d\xi_1.
\]

 \noi
 Then, we have the following bilinear estimate.
 See  Lemma 3.1  and  Corollary 3.2 in \cite{GV}.
 
\begin{lemma}
\label{LEM:bi-1}
Let $I^\frac{1}{2}$ and $I^\frac{1}{2}_-$
be as above \textup{(}with $\ta = \frac 12$\textup{)}.
Then, we have\footnote{We  use $a+$ (and $a-$) to denote $a + \eps$ (and $a - \eps$, respectively) for arbitrarily small $\eps \ll 1$,
where an implicit constant is allowed to depend on $\eps > 0$
(and it usually diverges as $\eps \to 0$).}

\begin{align*}
\big\| I^{\frac1{2}} I_-^{\frac1{2}} (u,v) \big\|_{ L^2_{x,t} (\R^2)} 
\les \| u\|_{X^{0,\frac1{2}+}} \| v\|_{X^{0,\frac1{2}+}}.
\end{align*}
\end{lemma}

The following two estimates are immediate corollary of Lemma \ref{LEM:bi-1}.

\begin{corollary}
\label{COR:bi}
\textup{(i)}
Let $N_1,N_2\geq 1$ be dyadic such that  $N_1 \gg N_2$. 
Then, we have 
\begin{align*}
\| \P_{N_1} u \P_{N_2} v \|_{L^2_{x, t}(\R^2)} 
\les \frac1{N_1}
 \| \P_{N_1}u \|_{X^{0,\frac1{2}+}}
 \| \P_{N_2}v \|_{X^{0,\frac1{2}+}}.
 \end{align*}

\noi
\textup{(ii)}
Let $m,n\in \Z$ such that  $|m+n|, |m-n|\ge 2$. Then, we have
\begin{align*}
\| \Pi_m u \Pi_n  v \|_{L^2_{x, t}(\R^2)} 
\les \frac1{\sqrt{|m+n||m-n|}}
 \|\Pi_m u\|_{X^{0,\frac1{2}+}}
 \| \Pi_n v\|_{X^{0,\frac1{2}+}}.
 \end{align*}

 \noi
\end{corollary}

\vspace{4mm}

In \cite{Tao}, Tao presented a proof of local well-posedness of
mKdV \eqref{mKdV} in $H^\frac{1}{4}(\R)$ based on the Fourier restriction norm method
by establishing the following trilinear estimate.

\begin{lemma}[Corollary 6.3 in \cite{Tao}]
\label{COR:tri-1}
Given small $\eps > 0$, there exists $C_\eps > 0$ such that 
\begin{align}
\label{tri}
\| \dx (u_1u_2u_3) \|_{X^{\frac14,-\frac12+2\eps}} 
\le C_\eps \prod_{j=1}^3 \| u_{j} \|_{X^{\frac14,\frac12+\eps}}.
\end{align}
\end{lemma}

In \cite{Tao}, the estimate \eqref{tri}
was stated with $-\frac 12 + \eps$ for the temporal regularity $b$
on the left-hand side.
It is, however, easy to see that the result
also holds true with $-\frac 12 + 2\eps$.

\section{Proof of Theorems \ref{THM:1} and \ref{THM:2}}
\label{SEC:tri}


\subsection{Trilinear estimate}

In view of the linear estimates in Lemma \ref{LEM:lin}, 
local well-posedness of \eqref{mKdV} (Theorem \ref{THM:1})
follows from a standard contraction argument once we prove the following trilinear estimate.

\begin{proposition}
\label{PROP:trilinear}
Let $s\ge \frac14$ and $2\le p < \infty$. Then, 
given small $\eps > 0$, there exists $C_\eps > 0$ such that 
\begin{align}
\big\| u_1 \cj{u}_2 \dx u_3 \big\|_{X^{s,-\frac12+2\eps}_{p}([0, T])} 
\le C_\eps  \prod_{j=1}^3\| u_{j} \|_{X^{s,\frac12+\eps}_{p}([0, T])}
\label{key}
\end{align}

\noi
for any $T > 0$.
\end{proposition}

\noi

We present the proof of Proposition \ref{PROP:trilinear} in the remaining part of this section.
By a standard reduction, it suffices to prove \eqref{key} without
the time restriction.
Noting that the resonance relation $\tau = \xi^3$ is invariant
under $(\tau, \xi) \mapsto (-\tau, -\xi)$, 
it suffices to prove 
\begin{align}
\big\| u_1 u_2 \dx u_3 \big\|_{X^{s,-\frac12+2\eps}_{p}} 
\les \prod_{j=1}^3\| u_{j} \|_{X^{s,\frac12+\eps}_{p}}.
\label{key0}
\end{align}

\noi
Furthermore, by the triangle inequality: $\jb{\xi} \les \jb{\xi_1}\jb{\xi_2}\jb{\xi_3}$
under $ \xi_1 + \xi_2 + \xi_3 +\xi = 0$, it suffices to prove~\eqref{key0}
for $s = \frac 14$.
Then, by duality, \eqref{key} follows once we prove
\begin{align}
\label{key1}
\bigg| \iint_{\R \times \R}
u_1 u_2 \dx u_3 \jb{\dx}^\frac{1}{4} v  dxdt \bigg| 
& = \bigg|  \intt_{\substack{\xi_1+\xi_2+\xi_3 + \xi =0 \\ \tau_1+\tau_2+\tau_3 + \tau =0} } 
\jb{\xi}^{\frac14} \xi_3
 \prod_{j=1}^3 \ft u_{j} (\xi_j,\tau_{j})\, \ft v(\xi, \tau) \bigg|\notag \\
& \les \prod_{j=1}^3 \| u_j\|_{X^{\frac14,\frac12+\eps}_{p}} \| v\|_{X^{0,\frac12-2\eps }_{p'}}.
\end{align}

\noi
In the following, 
we use 
$\xi_{\max}, \xi_{\med}, \xi_{\min}$ 
to denote the rearrangement of 
$\xi_1,\xi_2,\xi_3$ such that 
 $|\xi_{\max}|\ge  |\xi_{\med}|\ge |\xi_{\min}|$.
Under $\xi_1+\xi_2+\xi_3 + \xi =0$, 
we have $|\xi|\les |\xi_{\max}|$.
In the following, we apply dyadic decompositions
$|\xi_j|\sim N_j$ and $|\xi|\sim N$
for dyadic $N_j, N \geq 1$.
In this case, 
we also use the notation:
 $N_{\max}\sim |\xi_{\max}|$, 
  $N_{\med}\sim |\xi_{\med}|$, 
  and  $N_{\min}\sim |\xi_{\min}|$.

We prove Proposition  \ref{PROP:trilinear}
by separately considering the following four cases:

\begin{enumerate}

\item[(i)]
Trivial cases, 

\smallskip

\item[(ii)]
Non-resonant case:  $N_{\max} \gg N_{\med}$,

\smallskip

\item[(iii)]
Semi-resonant case: $N_{\max} \sim N_{\med}\gg  N_{\min}$,

\smallskip

\item[(iv)]
Resonant case: $N_{\max} \sim N_{\min}$.

\end{enumerate}

\noi

\noi
As we see below, the main difficulty appears in 
the resonant case (iv).
Before going into the details of the proof, 
we introduce a few more notations.
We use $\s$ and $\s_j$ to denote modulations given by 
\[
\s= \tau - \xi^3
\qquad \text{and} \qquad \s_j  = \tau_j - \xi_j^3
\]

\noi
for $j = 1,2,3$. We also  set
\[
\s_{\max} = \max \big(|\s|, |\s_1|, |\s_2|, |\s_3|\big).
\]

\noi
For conciseness of the presentation, 
we use the following (slightly abusive) short-hand notations:
\[u_N = \P_N u
\qquad\text{and}\qquad
u_n = \Pi_n u,\] 

\noi
where $\P_N$ is the Littlewood-Paley projector
and $\Pi_n$ is as in \eqref{Pi1}.
We only use the capitalized variables to denote dyadic numbers
and hence there is no confusion.

\begin{remark}\label{REM:pers}\rm
By slightly modifying the proof, we
can easily extend \eqref{key} to $1 \leq p < 2$.
Note that the proof in this case is easier than that of Proposition \ref{PROP:trilinear}
since $\l^p(\Z) \subset \l^2(\Z)$.
Furthermore, we can also establish 
\begin{align}
\big\| u_1 \cj u_2  \dx u_3 \big\|_{X^{s,-\frac12+2\eps}_{p}([0, T])} 
\le C_\eps  \min_{j = 1, 2, 3}\bigg( 
\| u_j \|_{X^{s,\frac12+\eps}_{p}([0, T])}
\prod_{\substack{k = 1\\k \ne j}}^3\| u_k \|_{X^{s,\frac12+\eps}([0, T])} \bigg)
\label{keyxx}
\end{align}

\noi
for  $s\ge \frac14$ and $1\le p < 2$. 
The tame estimate \eqref{keyxx} allows us to prove local well-posedness
of~\eqref{mKdV}
in $M^{2, p}_s(\R)$
for  $s\ge \frac14$ and $1\le p < 2$,
where 
the local existence time depends only
on the $H^s$-norm of initial data.
In particular, this allows us to prove global well-posedness
of~\eqref{mKdV}
in $M^{2, p}_s(\R)$
for  $s\ge \frac14$ and $1\le p < 2$.
See Appendix of \cite{OW2} for such an argument.
Since the required modification is straightforward, 
we omit details.

\end{remark}

\subsection{Trivial cases}
We first consider two trivial cases:
\begin{align}
\text{(i)} \quad |\xi_{\max}|\les 1
\qquad \text{and} \qquad \text{(ii)} \quad  \jb{\s_{\max}} \gg \jb{\xi_{\max}}^{10}.
\label{trivial}
\end{align}

\noi
(i) 
Suppose $|\xi_{\max}|\les 1$.
In this case, we have $|\xi | \les 1$.
Then, by H\"older's inequality,  Bernstein's inequality \eqref{bern},  \eqref{Xsb}, \eqref{embed1}
followed by \eqref{embed1a} and \eqref{embed2}, 
we have
\begin{align*}
 \sum_{N_{\max}, N  \les 1} 
 & N^{\frac14}\bigg|\int_{\R\times\R} 
 u_{N_1} u_{N_2} \dx u_{N_3}   v_N dxdt \bigg|\\
& \les   \sum_{N_{\max}, N \les 1}   
\|u_{N_1}\|_{L^2_{t, x}} \|u_{N_2}\|_{L^\infty_{x, t}} \|{u_{N_3}}\|_{L^\infty_{x, t}}  \|v_N\|_{L^2_{x, t}}\\
& \les   \sum_{N_{\max}, N  \les 1}  
 \|u_{N_1}\|_{L^2_{x, t}} \|u_{N_2}\|_{L^\infty_tL^2_x} \|{u_{N_3}}\|_{L^\infty_tL^2_x}  \|v_N\|_{L^2_{x, t}}\\
& \les  \sum_{N_{\max}, N  \les 1} \bigg( \prod_{j = 1}^3  \|u_{N_j}\|_{X^{0,\frac12+}} \bigg)
 \|v_N\|_{X^{0,\frac12-}}\\
& \les  \sum_{N_{\max}, N  \les 1} 
\bigg( \prod_{j = 1}^3   \|u_{N_j}\|_{X^{0,\frac12+}_{p}} \bigg)  \|v_N\|_{X^{0,\frac12-}_{p'}}.
\end{align*}

\noi
By summing over dyadic blocks $N_1, N_2, N_3, N \les 1$,
we obtain \eqref{key1}.

\smallskip

\noi
(ii) 
Next, we suppose  $ \jb{\s_{\max}} \gg \jb{\xi_{\max}}^{10}$.
In the following, we consider the case $\jb{\s_1} = \jb{\s_{\max}}$,
The other cases follow from a similar argument.
By H\"older's and  Bernstein's inequalities, the definition \eqref{Xsb}, 
and \eqref{embed1}, we have
\begin{align*}
\sum_{\substack{N_1, N_2,N_3 N \geq 1\\\text{dyadic}}} & N^{\frac14} 
\bigg|\int_{\R\times \R} u_{N_1} u_{N_2} \dx u_{N_3}   v_N dxdt \bigg|\\
& \les  \sum_{N_1, N_2,N_3, N\geq 1}   N_{\max}^{\frac54} \|u_{N_1}\|_{L^2_{x, t}} \|u_{N_2}\|_{L^\infty_{x, t}} \|u_{N_3}\|_{L^\infty_{x, t}} \|v_N\|_{L^2_{x, t}}\\
& \les  \sum_{N_1, N_2,N_3, N \geq 1}   N_{\max}^{\frac94} \|u_{N_1}\|_{L^2_{x, t}} \|u_{N_2}\|_{L^\infty_{t}L^2_x} \|u_{N_3}\|_{L^\infty_{t}L^2_x} \|v_N\|_{L^2_{x, t}}\\
& \les  \sum_{N_1, N_2,N_3, N  \geq 1}   N_{\max}^{\frac94} \jb{\s_1}^{-\frac12-} 
\bigg(\prod_{j = 1}^3 \|u_{N_j}\|_{X^{0,\frac12+}}\bigg) 
 \|v_N\|_{X^{0,\frac12-}}
\intertext{By applying the lower bound \eqref{trivial}
together with \eqref{embed1a} and \eqref{embed2}, }
& \les   \sum_{N_1, N_2,N_3, N \geq 1}   N_{\max}^{\frac94} N_{\max}^{-\frac{7}{2} - \frac{3}{p} -} 
\bigg(\prod_{j = 1}^3 \|u_{N_j}\|_{X^{0,\frac12+}_p}\bigg) 
\|v_N\|_{X^{0,\frac12-}_{p'}}\\
& \les   \sum_{N_1, N_2,N_3, N \geq 1}  
\bigg(\prod_{j = 1}^3 N_j^{0-} \|u_{N_j}\|_{X^{0,\frac12+}_p}\bigg) 
N^{0-} \|v_N\|_{X^{0,\frac12-}_{p'}}.
\end{align*}

\noi
By summing over dyadic blocks $N_1, N_2, N_3, N \geq 1$,
we obtain \eqref{key1}.

Therefore, we assume that 
\begin{align}
|\xi_{\max}| \gg 1 \qquad \text{and}\qquad \jb{\s_{\max}} \les \jb{\xi_{\max}}^{10}
\label{modbd}
\end{align}

\noi
in the following.

\begin{remark}\rm
\label{REM:strategy}
In the  arguments above, 
we first established bounds in terms of the  standard $X^{s,b}$-norms
and then applied \eqref{embed1a} and \eqref{embed2}
to replaced it by the $X^{s, b}_p$-norms.
More precisely, we used
\begin{align}
\label{holder1}
\|u_N\|_{X^{s,b}} 
 \les  \max \big( N^{\frac12 -\frac1p}, 1 \big) \|u_{N}\|_{X^{s,b}_{p}}
\end{align}

\noi
and
\begin{align}
\label{holder2}
\sum_{\substack{N \geq 1\\\text{dyadic}}} N^{-\eps} \|u_{N}\|_{X^{s,b}_{p}} \les \|u\|_{X^{s,b}_{p}} 
\end{align}

\noi
for any $\eps > 0$.
We use the same strategy in the following.
\end{remark}

\subsection{Non-resonant case: $N_{\max} \gg N_{\med} \ge N_{\min}$}
\label{SUBSEC:non}

Without loss of generality,\footnote{Since the derivative falls on the third factor
on the left-hand side of \eqref{key1}, 
there is no symmetry among frequencies $\xi_1, \xi_2$, and $\xi_3$.
However, we simply bound this derivative by the largest frequency in the following
and thus we may pretend that there is symmetry among frequencies.
The same comment applies in the following.} suppose that $N_1 \gg N_2 \geq N_3$.
The other cases can be treated by a similar consideration.
In this case, we have  $N \sim N_1$. 
Then, by  Corollary \ref{COR:bi} and \eqref{modbd}, we have
\begin{align}
 \sum_{N_1\sim N\gg N_2}
&  N^{\frac14} \bigg| \int_{\R\times \R}u_{N_1} u_{N_2} \dx u_{N_3}  v_N dxdt \bigg| \notag \\
& \les  \sum_{N_1\sim N\gg N_2}   N^{\frac54} \|u_{N_1} u_{N_2}\|_{L^2_{x, t}} \|u_{N_3} v_N\|_{L^2_{x, t}}\notag \\
& \les \sum_{N_1\sim N\gg N_2}  N^{-\frac34} 
\bigg(\prod_{j = 1}^3  \| u_{N_j} \|_{X^{0,\frac12+}} \bigg)  \| v_{N} \|_{X^{0,\frac12+}} \notag \\
& \les \sum_{N_1\sim N\gg N_2}  N^{-\frac34+}
\bigg(\prod_{j = 1}^3  \| u_{N_j} \|_{X^{0,\frac12+}} \bigg)  \| v_{N} \|_{X^{0,\frac12-}} \notag \\
\intertext{By applying \eqref{holder1} and \eqref{holder2}, }
& \les \sum_{N_1\sim N\gg N_2}  N^{-\frac12-\frac1p+} (N_2N_3)^{\frac14-\frac1p} 
\bigg(\prod_{j = 1}^3  \| u_{N_j} \|_{X^{\frac14,\frac12+}_{p}}\bigg)
\| v_{N} \|_{X^{0,\frac12-}_{p'}}   \notag \\
& \les \sum_{N_1\sim N\gg N_2}  N^{-\frac1p+} 
\bigg(\prod_{j = 1}^3  \| u_{N_j} \|_{X^{\frac14,\frac12+}_{p}}\bigg)
 \| v_{N} \|_{X^{0,\frac12-}_{p'}} \notag \\
& \les \prod_{j=1}^3 \|u_j \|_{X^{\frac14,\frac12+}_{p}}  \|v \|_{X^{0,\frac12-}_{p'}} , 
\label{X1}
\end{align}

\noi
provided $p< \infty$.

\subsection{Semi-resonant case: $N_{\max} \sim N_{\med} \gg N_{\min}$}
\label{SUBSEC:half}
We proceed as in the non-resonant case.
The frequency separation allows us to use the  bilinear estimate (Corollary \ref{COR:bi}) twice,
gaining two derivatives.
Without loss of generality, suppose  that $N_1\sim N_2 \gg N_3$.
The other cases can be treated by a similar consideration.
We distinguish two cases according to the relation between $N$ and $N_{\max}$.

First, suppose that $N \ll N_{\max}$. Then,
by Corollary \ref{COR:bi}, we have
\begin{align*}
\sum_{N_1\sim N_2 \gg N_3,   N}  & N^{\frac14} 
 \bigg|\int_{\R\times \R} u_{N_1} u_{N_2} \dx u_{N_3}  v_N dxdt \bigg|\\
& \les  \sum_{N_1\sim N_2 \gg N_3,   N}   N_{\max}^{\frac54} \|u_{N_1} u_{N_3}\|_{L^2_{x, t}}
 \|u_{N_2} v_N\|_{L^2_{x, t}}\\
& \les \sum_{N_1\sim N_2 \gg N_3,   N}  N_{\max}^{-\frac34} 
\bigg( \prod_{j = 1}^3  \| u_{N_j} \|_{X^{0,\frac12+}} \bigg) \| v_{N} \|_{X^{0,\frac12+}}.
\end{align*}

\noi
The rest follows as in \eqref{X1}.

Next, consider the case $ N\sim  N_{\max}$.
In this case, we have  $|\xi_1+ \xi_2 + \xi|  = |\xi_3| \ll N\sim N_{\max}$. 
Hence, we must have $\xi_1\xi_2 < 0$, $\xi_1\xi < 0$, or $\xi_2 \xi < 0$.
Without loss of generality, suppose that $\xi_1\xi_2 <0$.
(The  proofs for the other cases are similar.) 
We then have 
$|\xi - \xi_3||\xi + \xi_3| \sim N_{\max}^2$
and 
\[|\xi_1-\xi_2||\xi_1+\xi_2|
= |\xi_1-\xi_2||\xi+\xi_3| \sim N_{\max}^2.\]

\noi
Hence, by Lemma \ref{LEM:bi-1}, 
we have
\begin{align}
 \sum_{N_1\sim N_2 \sim N \gg N_3} & N^{\frac14} 
\bigg| \int_{\R\times \R} u_{N_1} u_{N_2} \dx u_{N_3}  v_N dxdt \bigg|\notag \\
& \les  \sum_{N_1\sim N_2 \sim N \gg N_3}   N_{\max}^{\frac54} \|u_{N_1} u_{N_2}\|_{L^2} \|u_{N_3} v_N\|_{L^2}\notag \\
& \les \sum_{N_1\sim N_2 \sim N \gg N_3}  N_{\max}^{-\frac34}
\bigg( \prod_{j = 1}^3   \| u_{N_j} \|_{X^{0,\frac12+}} \bigg)  \| v_{N} \|_{X^{0,\frac12+}}.
\label{X2}
\end{align}

\noi
The rest follows as in \eqref{X1}.

\subsection{Resonant case}
\label{SUBSEC:res}
In this case,
we have $N_1\sim N_2 \sim N_3$.
Without loss of  generality,
we may further assume that $N_1\sim N$, since, otherwise, i.e.~$N_1 \gg N$, 
the proof can be reduced to \eqref{X2} with the roles of $N$ and $N_3$ switched.

Hence, we assume that $N_1 \sim N_2 \sim N_3 \sim N$ in the following.
This case requires more careful analysis 
and we need to use the unit-cube decomposition:
\[ u = \sum_{n \in \Z} u_n = \sum_{ n \in \Z} \Pi_n u. \]

\noi
Given $n \in \Z$, we set $I_n = [n, n+1)$.

\smallskip

\noi
$\bullet$ {\bf Case 1:}
We first consider the case  $|\xi_i - \xi_j | \ge |\xi_i + \xi_j|$
for some pair   $(i,j)$.
\\
\indent
Without loss of  generality,
 we assume $(i,j) = (1,2)$.
In the next two subcases, 
we treat the  case $|\xi_1 + \xi_2| \les 1$.

\smallskip

\noi
\underline{\bf Subcase 1.1:} $|\xi_1 + \xi_2| \les 1$
and $\min( |\xi_1 -\xi_3|, |\xi_1 +\xi_3|) \les 1$.

We only consider the case where $|\xi_1 - \xi_3|\les 1$, 
since  the proof for the case  $|\xi_1 + \xi_3|\les 1$ is similar.
Suppose that $\xi_1 \in I_n = [n, n+1)$.
Then, we have
\[\xi_2 = -n + O(1), \qquad \xi_3 = n + O(1), \qquad
\text{and}\qquad\xi = -n + O(1).\]

\noi
Hence, we need to estimate the following expression:
\begin{align*}
\sum_{n\in \Z} \sum_{j, k, \l = O(1)} \jb{n}^{\frac14} 
\bigg|\int_{\R\times \R} u_{n} u_{-n+j} \dx u_{n+k}  v_{-n+\l} dxdt \bigg|.
\end{align*}

\noi
For simplicity of the presentation, 
 we only consider the ``diagonal'' case, i.e.
$j = k = \l = 0$ in the following.
By  Lemma  \ref{COR:tri-1}
and H\"older's inequality in $n$,  
we have
\begin{align*}
\sum_{n\in \Z}  \jb{n}^{\frac14} 
\bigg|\int_{\R\times \R} u_{n} u_{-n} \dx u_{n}  v_{-n} dxdt \bigg|
& \les \sum_{n\in \Z}   \jb{n}^{\frac14} \| u_{n} u_{-n} \dx u_n \|_{X^{0,-\frac12+2\eps}}
 \| v_{-n}\|_{X^{0,\frac12-2\eps}}\\
& \les \sum_{n\in \Z}  \|u_n\|_{X^{\frac14,\frac12+\eps}}^3 \| v_{n}\|_{X^{0,\frac12-2\eps}}\\
& \les  \|u\|_{X^{\frac14,\frac12+\eps}_{3p}}^3  \| v\|_{X^{0,-\frac12-2\eps}_{p'}} \\
& \les  \|u\|_{X^{\frac14,\frac12+\eps}_{p}}^3  \| v\|_{X^{0,-\frac12-2\eps}_{p'}}
\end{align*}

\noi
for sufficiently small $\eps > 0$.
This is the only case where we need to be precise about  the temporal regularities.

\smallskip

\noi
\underline{\bf Subcase 1.2:} $|\xi_1 + \xi_2| \les 1$ and $|\xi_1 \pm \xi_3| \gg 1$.
\\
\indent
Suppose that $\xi_1 \in I_n$
and $\xi_3 \in I_m$.
Then, we have
$\xi_2 \in I_{-n + j}$ and $\xi \in I_{-m + k}$
for $j, k = O(1)$.
As in Subcase 1.1, we only estimate the contribution from  $j = k = 0$.
Without loss of generality, we assume that $|m+n| \geq |m-n|$.
By Corollary \ref{COR:bi} with $|m \pm n| \gg 1$, 
 $|m|\sim |n| \gg 1$, and \eqref{modbd}, 
 we have
\begin{align*}
\sum_{\substack{m, n \in \Z\\m\ne n}} 
& |m|^{\frac14} \bigg| \int_{\R\times\R} u_{n} u_{-n} \dx u_m  v_{-m} dxdt \bigg|
 \leq \sum_{\substack{m, n \in \Z\\m\ne n}}
  |m|^{\frac54} \| u_{m} u_{n}  \|_{L^2} \| u_{-n} v_{-m}\|_{L^2}\\
& \les \sum_{\substack{m, n \in \Z\\m\ne n}}
  \frac{|m|^{\frac54}}{|m-n||m+n|} \|u_n\|_{X^{0,\frac12+}} \|u_m\|_{X^{0,\frac12+}} 
  \|u_{-n}\|_{X^{0,\frac12+}}
  \| v_{-m}\|_{X^{0,\frac12+}}\\
& \les \sum_{\substack{m, n \in \Z\\m\ne n}}
  \frac{|n|^{\frac12+}}{|m-n||n+m|} \|u_n\|_{X^{\frac14,\frac12+}} \|u_m\|_{X^{\frac14,\frac12+}}
   \|u_{-n}\|_{X^{\frac14,\frac12+}}
   \| v_{m}\|_{X^{0,\frac12-}}\\
& \les  \sum_{\substack{m, n \in \Z\\m\ne n}}
 \frac{ 1}{|m-n|\jb{n}^{\frac 12 -}} \|u_n\|_{X^{\frac14,\frac12+}} \|u_m\|_{X^{\frac14,\frac12+}} 
 \|u_{-n}\|_{X^{\frac14,\frac12+}}
 \| v_{m}\|_{X^{0,\frac12-}}
\intertext{By applying  \eqref{hlp2} and \eqref{Xsb}, }
& \les  \Big\|  \|u_n\|_{X^{\frac14,\frac12+}}
 \|u_{-n}\|_{X^{\frac14,\frac12+}} \Big\|_{\l^{\frac{p}{2}}_n} 
  \Big\| \|u_m\|_{X^{\frac14,\frac12+}} \| v_{m}\|_{X^{0,\frac12-}} \Big\|_{\l^{\frac{p}{p-2}}_m}\\
& \les  \|u\|_{X^{\frac14,\frac12+}_{p}}^3  \| v\|_{X^{0,\frac12+}_{p'}}.
\end{align*}

In the next three subcases, we treat the case $|\xi_1 + \xi_2| \gg 1$.

\smallskip

\noi
\underline{\bf Subase 1.3:} $|\xi_1+\xi_2| \gg 1$ and $|\xi_i -  \xi_j| \les 1$ for some $(i,j) \neq (1,2)$.
\\
\indent
Without loss of generality, we may assume $(i,j) = (1,3)$.
Suppose that $|\xi - \xi_3|\les 1$.
Then, we need to show
\[
\sum_{n\in \Z} \jb{n}^{\frac14} 
\bigg|\int_{\R\times\R} u_{n} u_{-3n} \dx u_n  v_{n} dxdt \bigg|
\les \|u\|_{X^{\frac14,\frac12+}_{p}}^3  \| v\|_{X^{0,\frac12+}_{p'}},
\]

\noi
which can be easily obtained by repeating the argument in Subcase 1.1.
Hence, we assume that 
 $|\xi-\xi_3| \gg 1$ in the following.

Suppose that $\xi_1 \in I_n$ and $\xi \in I_m$.
Then, we have 
$\xi_3 \in I_{n+ j}$ and $\xi_2 \in I_{-m-2n-k}$
for $j, k = O(1)$.
 As above, we only estimate the  contribution from $j = k = 0$.
By the triangle inequality, we have 
 $\max(|\xi - \xi_3|, |\xi+\xi_3|) \ges |m| \sim |n| \gg 1$.
In the following, we only consider the case
$ |\xi - \xi_3| \sim |m|$ since the other case follows in a similar manner.
Moreover, since $|\xi_1 - \xi_2 | \geq |\xi_1 + \xi_2|$, 
we conclude that 
$|m + 3n| \sim |m|$.
Hence, 
by Corollary~\ref{COR:bi}, 
 \eqref{modbd}, and 
 \eqref{hlp2}, we have
\begin{align*}
\sum_{\substack{m, n \in \Z\\m\ne n}} 
& |m|^{\frac14} \bigg|\int_{\R\times\R} u_{n} u_{-m-2n} \dx u_n  v_{m} dxdt\bigg| \\
& \les \sum_{\substack{m, n \in \Z\\m\ne n}} 
 |m|^{\frac54} \| u_{n} u_{-m-2n}  \|_{L^2} \| u_{n} v_{m}\|_{L^2}\\
& \les \sum_{\substack{m, n \in \Z\\m\ne n}} 
  \frac{|m|^{\frac12+}}{|m+n|\sqrt{|m-n||m+3n|}} 
\| u_{n}\|_{X^{\frac14,\frac12+}}^2 \|u_{-m-2n}  \|_{X^{\frac14,\frac12+}}
 \|v_{m}\|_{X^{0,\frac12-}}\\
& \les
 \| u\|_{X^{\frac14,\frac12+}_\infty}^2
 \sum_{\substack{m, n \in \Z\\m\neq n}}  \frac{1}{|m + (m - 2n) |\jb{m}^{\frac 12-}}
 \|u_{m}  \|_{X^{\frac14,\frac12+}}  \|v_{m-2n}\|_{X^{0,\frac12-}}\\
& \les  \|u\|_{X^{\frac14,\frac12+}_{p}}^3  \| v\|_{X^{0,\frac12+}_{p'}}.
\end{align*}

\smallskip

\noi
\underline{\bf Subcase 1.4:} $|\xi_1+\xi_2| \gg 1$ and $|\xi_i + \xi_j| \les 1$ for some $(i,j) \neq (1,2)$.
\\
\indent
We can proceed as in Subcase 1.3 above
and thus we omit details.


\smallskip

\noi
\underline{\bf Subcase 1.5:} $|\xi_1+\xi_2| \gg 1$ and $|\xi_i \pm \xi_j| \gg 1$ for all $(i,j) \neq (1,2)$.
\\
\indent
By assumption, we have  $|\xi_1 - \xi_2|\ge |\xi_1 + \xi_2| $
and hence
we have  $|\xi_i \pm \xi_j| \gg1$ for all $i\neq j$.
Recalling that 
\[
\s + \s_1 + \s_2 + \s_3 = 3(\xi_1+\xi_2)(\xi_2+ \xi_3)(\xi_1 + \xi_3)
\]

\noi
under $\xi_1 + \xi_2 + \xi_3 + \xi = 0$
and $\tau_1+\tau_2+\tau_3+\tau= 0$, 
we have
\begin{align}
\label{modu}
\jb{\s_{\max}} \ges |\xi_1+\xi_2||\xi_2+ \xi_3||\xi_1 + \xi_3|.
\end{align}

\noi
Without loss of generality, we assume that $\jb{\s_1} = \jb{\s_{\max}}$.
By  Bernstein's inequalities,~\eqref{modu},  and Corollary \ref{COR:bi} with \eqref{modbd}, 
we have
\begin{align}
 \sum_{n_1+n_2+n_3+n = O(1)} 
&  |n|^{\frac14} \bigg|\int_{\R\times\R} u_{n_1} u_{n_2} \dx u_{n_3} v_{n} dxdt\bigg| \notag \\
& \les \sum_{n_1+n_2+n_3+n =O(1)}  |n|^{\frac54} \| u_{n_1}\|_{L^2_{x, t}} \| u_{n_2} \|_{L^\infty_{t}L^2_x} \| u_{n_3} v_{n}\|_{L^2_{x, t}}\notag \\
& \les \sum_{n_1+n_2+n_3+n = O(1)}  \frac{|n|^{\frac54+}}{\sqrt{|n_1+n_2||n_1+n_3||n_2 + n_3|} 
\sqrt{|n_3-n||n_3+n|}} \notag \\
& \hphantom{XXXXXX}
\times \bigg(\prod_{j = 1}^3 \|u_{n_j}\|_{X^{0,\frac12+}}\bigg)
\| v_{n}\|_{X^{0,\frac12-}}.
\label{X3}
\end{align}

\noi
By the triangle inequality, we have
$\max(|n_3-n|, |n_3+n|) \geq |n|$
and 
\[
\max(|n_1+n_3|, |n_2 + n_3|)
\geq |n_1 - n_2| \ges |n|.\]

\noi
In the following, we only consider the case
 $|n_1+n_3| \sim |n_3-n| \ges |n|$.
Then, we have
\begin{align*}
\text{LHS of \eqref{X3}}
& \les \sum_{n_1+n_2+n_3+n = O(1) }  
\frac{|n|^{-\frac12+}}{\sqrt{|n_1+n_2||n_2 + n_3||n_3+n|}} \\
& \hphantom{X}
\times \bigg(\prod_{j = 1}^3 \|u_{n_j}\|_{X^{0,\frac12+}}\bigg)
\| v_{n}\|_{X^{0,\frac12-}}\\
& \les  \sup_{n_3, n}
\bigg( \sum_{n_2}  \frac{|n_2|^{-\frac12+}}{\sqrt{|n_2 + n_3|}} \|u_{-n_2-n_3-n}\|_{X^{\frac14,\frac12+}} \|u_{n_2}\|_{X^{\frac14,\frac12+}} \bigg)  \\
& \hphantom{X}
 \times \sum_{n_3,n}  \frac{1}{|n_3+n|\jb{n}^{0+}}  \|u_{n_3}\|_{X^{\frac14,\frac12+}} \| v_{n}\|_{X^{0,\frac12-}}\\
\intertext{By  applying H\"older's inequality in $n_2$ and \eqref{hlp2},}
& \les  \|u\|_{X^{\frac14,\frac12+}_{p}}^3  \| v\|_{X^{0,\frac12+}_{p'}},
\end{align*}

\noi
provided that $p < \infty$.

\smallskip

\noi
$\bullet$ {\bf Case 2:}
$|\xi_i - \xi_j | \le |\xi_i + \xi_j|$ for all $i,j$.

In this case, all $\xi_j$'s for $j =1,2,3$ have the same sign.
Thus, we have
$|\xi \pm \xi_j| \ges |\xi_{\max}|$ for $j=1,2,3$.
Moreover,  from \eqref{modu}, 
we have  
\begin{align}
\jb{\s_{\max}}\ges |\xi_{\max}|^3.
\label{X4}
\end{align}

We first consider the case
$\s_j = \s_{\max}$ for some $j = 1, 2, 3$.
Without loss of generality, we  assume that  $\s_2 = \s_{\max}$.
By H\"older's and Bernstein's inequalities, 
\eqref{Xsb}, \eqref{X4}, and Lemma~\ref{LEM:bi-1}
with $|\xi \pm \xi_3| \ges |\xi_{\max}|$,   
we have 
\begin{align*}
&\sum_{N_{\max} \sim N_{\min}\sim N} N^{\frac54} 
\bigg|\int_{\R\times\R} u_{N_1} u_{N_2} u_{N_3} v_{N}  dxdt  \bigg| \\
& \les \sum_{N_{\max} \sim N_{\min}\sim N} N^{\frac74} \| u_{N_1}\|_{L^\infty_{t}L^2_x}  \| u_{N_2}\|_{L^2_{x, t}} \| u_{N_3} v_{N}\|_{L^2_{x, t}}   \\
& \les \sum_{N_{\max} \sim N_{\min}\sim N } N^{\frac14-} \| u_{N_1}\|_{X^{0,\frac12+}} 
 \| u_{N_2}\|_{X^{0,\frac12+}} \| u_{N_3} v_{N}\|_{L^2_{x, t}}   \\
& \les \sum_{N_{\max} \sim N_{\min}\sim N} N^{-\frac34+} 
\bigg(\prod_{j = 1}^3 \| u_{N_j}\|_{X^{0,\frac12+}} \bigg) \| v_{N}\|_{X^{0,\frac12-}}.
\end{align*}

\noi
Then, the rest follows as in \eqref{X1}.

In the following, we assume that 
 $\s = \s_{\max}$. The proof for this case is more involved 
 and 
thus we split it into several subcases.

\smallskip

\noi
\underline{\bf Subcase 2.1:}
$\s = \s_{\max}$ and $|\xi_i - \xi_j| \les1$ for some $i\ne j$.
 
Without loss of generality, we may assume $|\xi_1 - \xi_2| \les 1$.
Suppose that $\xi_1 \in I_n$ and $\xi_3 \in I_m$.
Then, we have 
$\xi_2 \in I_{n+ j}$ and $\xi_3 \in I_{-m-2n-k}$
for $j, k = O(1)$.
In the following, 
we only estimate the  contribution from  $j = k = 0$:
\begin{align}
\label{diag1}
&\sum_{\substack{m, n \in \Z\\|n|\sim |m|}} 
|n|^{\frac54} \bigg|\int_{\R\times \R} u_{n} u_{n} u_{m} v_{-m-2n}  dxdt  \bigg|. 
\end{align}

We first consider the case $|\xi_1 - \xi_3| \les1$.
In this case, we can further reduce \eqref{diag1}
 to the following diagonal case:
\begin{align}
\sum_{n\in \Z} |n|^{\frac54} \bigg|\int_{\R\times \R} u_{n} u_{n} u_{n} v_{-3n}  dxdt \bigg|.
\label{diag2}
\end{align}

\noi
By H\"older's inequality,  Bernstein's inequality \eqref{bern} 
and  \eqref{X4}, we have
\begin{align*}
\eqref{diag2}
& \les \sum_{n\in \Z} |n|^{\frac54}  \| u_n\|_{L^\infty_{t}L^2_x}^2  \| u_n\|_{L^2_{x, t}}  \| v_{-3n}\|_{L^2_{x, t}} \\
& \les \sum_{n\in \Z} |n|^{-\frac14+}  \| u_{n}\|_{X^{0,\frac12+}}^3
\| v_{3n}\|_{X^{0,\frac12-}} .
\end{align*}

\noi
Then, the rest follows from H\"older's inequality in $n$.

Next, we consider the case 
$|\xi_1 - \xi_3| \gg1$. 
In this case, we have $|m + n | \ge |m - n| \gg 1$.
By H\"older's and Bernstein's inequalities, 
\eqref{X4}, Corollary \ref{COR:bi}, we have
\begin{align*}
\eqref{diag1}
& \les \sum_{m, n \in \Z}
 |n|^{\frac54} \| u_{m} u_{n}\|_{L^2_{x, t}} \|u_{n}\|_{L^\infty_{t}L^2_{x}} \|v_{-m-2n}\|_{L^2_{x, t}}   \\
& \les\sum_{m, n \in \Z}
\frac{ |n|^{-\frac14+} }{\sqrt{|m-n||m+n|}} \| u_{n}\|_{X^{0,\frac12+}}^2  \| u_{m}\|_{X^{0,\frac12+}}
 \|v_{-m-2n}\|_{X^{0,\frac12-}}   \\
& \les \|u\|_{X^{\frac14,\frac12+}_{\infty}}^2 
\sum_{m,  n\in \Z} \frac{1}{|n + (-m - 2n)|^\frac{1}{2} \jb{n}^{1-}} 
  \|u_{n}\|_{X^{\frac14,\frac12+}} \|v_{-m-2n}\|_{X^{0,\frac12-}}\\
  & \sim \|u\|_{X^{\frac14,\frac12+}_{\infty}}^2 
\sum_{m,  n\in \Z} \frac{1}{|n + (-m - 2n)| \jb{n}^{\frac 12-}} 
  \|u_{n}\|_{X^{\frac14,\frac12+}} \|v_{-m-2n}\|_{X^{0,\frac12-}}.
\end{align*}

\noi
Then, the rest follows from \eqref{hlp2}.

\smallskip

\noi
\underline{\bf Subcase 2.2:}
$\s = \s_{\max}$ and $|\xi_i - \xi_j| \gg1$ for all $i\ne j$.

Since all $\xi_j$'s have the same sign, we have $|\xi_i + \xi_j|\sim |\xi_i| \sim |\xi_{\max}|$. 
Then, by H\"older's and Bernstein's inequalities, 
\eqref{X4}, and 
Corollary \ref{COR:bi}  with  $|n_1\pm n_2| \gg 1$,
we have 
\begin{align}
& \sum_{n_1+n_2+n_3+n =O(1)} 
 |n|^{\frac54}  \bigg|\int_{\R\times\R}u_{n_1} u_{n_2} u_{n_3} v_{n}  dxdt  \bigg|  \notag \\
& \hphantom{XXX}
\les \sum_{n_1+n_2+n_3+n = O(1)} \jb{n}^{\frac54} \| u_{n_1} u_{n_2}\|_{L^2_{x, t}} \|u_{n_3}\|_{L^\infty_{t}L^2_{x}} \|v_{n}\|_{L^2_{x, t}}   \notag \\
& \hphantom{XXX}
 \les \sum_{n_1+n_2+n_3+n =O(1)} \jb{n}^{-\frac14+} \| u_{n_1} u_{n_2}\|_{L^2_{x, t}} 
\|u_{n_3}\|_{X^{0,\frac12+}} \|v_{n}\|_{X^{0,\frac12-}}   \notag\\
& \hphantom{XXX}
 \les \sum_{n_1+n_2+n_3+n = O(1) } \frac{\jb{n}^{-1+} }{\sqrt{|n_1-n_2||n_1+n_2|}} 
\bigg(\prod_{j = 1}^3 \| u_{n_j}\|_{X^{\frac14,\frac12+}} \bigg) \|v_{n}\|_{X^{0,\frac12-}}.
\label{X5}
\end{align}

\noi
By noting  $|n_1+n_2| \sim |n_3+n| \sim |n| \sim |n_1|$
and applying H\"older's inequality in $n_1$ and~\eqref{hlp2}, 
we have
\begin{align*}
\text{LHS of \eqref{X5}}
& \les \sum_{n_1+n_2+n_3+n =O(1)} \frac{\jb{n_1}^{-\frac12+}  \jb{n}^{0-}}{\sqrt{|n_1-n_2|}|n_3+n|} 
\bigg(\prod_{j = 1}^3 \| u_{n_j}\|_{X^{\frac14,\frac12+}}\bigg) \|v_{n}\|_{X^{0,\frac12-}_{2,2}}  \\
& \les\sup_{n_3,n} \bigg( \sum_{n_1}\frac{ \jb{n_1}^{-\frac12+} }{\sqrt{\jb{2n_1+n+n_3}}} \| u_{n_1}\|_{X^{\frac14,\frac12+}}  \| u_{-n_1-n_2-n}\|_{X^{\frac14,\frac12+}} \bigg) \\
& \hphantom{XX}
\times \left(  \sum_{n_3,n}  \frac{\jb{n}^{0-}}{|n_3+n|} \|u_{n_3}\|_{X^{\frac14,\frac12+}} \|v_{n}\|_{X^{0,\frac12-}}  \right) \\
& \les  \|u\|_{X^{\frac14,\frac12+}_{p}}^3  \| v\|_{X^{0,\frac12+}_{p'}},
\end{align*}

\noi
provided that $p < \infty$.

This completes the proof of Proposition \ref{PROP:trilinear}
and hence the proof of Theorem \ref{THM:1}.

\subsection{Persistence of regularity}
\label{SUBSEC:GWP}
We conclude this section 
by presenting the proof of global well-posedness (Theorem \ref{THM:2}).
When $\frac 14 \leq  s < 1- \frac{1}{p}$, 
global well-posedness immediately follows
from the local well-posedness in Theorem \ref{THM:1}
together with the global-in-time a priori bound \eqref{bd1}
in Proposition \ref{PROP:global}.
In the following, 
we briefly discuss the situation for $s \geq 1 - \frac 1p$.
In this case, 
the proof is based on combining 
the global-in-time a priori bound~\eqref{bd1}
in Proposition \ref{PROP:global}
on the $M^{2, p}_\frac{1}{4}$-norms of solutions
and a persistence-of-regularity argument.

With the notations from the previous subsections, 
we have $|\xi| \les |\xi_{\max}|$.
Hence, by slightly modifying the proof of Proposition \ref{PROP:trilinear},
we obtain 
\begin{align}
\big\| |u|^2 \dx u \big\|_{X^{s,-\frac12+2\eps}_{p}([0, T])} 
\le C_\eps  \| u \|_{X^{\frac{1}{4},\frac12+\eps}_{p}([0, T])}^2
\| u \|_{X^{s,\frac12+\eps}_{p}([0, T])}
\label{X6}
\end{align}

\noi
for any $s \geq \frac{1}{4}$, any $T > 0$,
and for small $\eps > 0$.

Let  $u_0 \in M^{2, p}_s(\R)$
for some  $s\geq \frac 14$ and $2\leq p < \infty$.
Since $u_0 \in M^{2, p}_\frac{1}{4}(\R)$, 
there exists a unique global solution
$u \in C(\R; M^{2, p}_\frac{1}{4}(\R))$  to~\eqref{mKdV}
with $u|_{t = 0} = u_0$.
We need to check that $u$ indeed lies in the correct space 
$ C(\R; M^{2, p}_s(\R))$.
In view of  the global-in-time a priori bound \eqref{bd1}, 
there exists small local existence time 
\begin{align}
\dl \sim (1+ \|u_0\|_{M^{2, p}_\frac{1}{4}})^\ta > 0
\label{X6a}
\end{align}

\noi
for some $\ta > 0$
such that 
a standard contraction argument
in $X^{\frac{1}{4},\frac12+\eps}_{p}(I)$
can be applied on any interval $I$of length $\dl$.
Moreover, with $I = [t_0, t_0 + \dl]$, we have
\begin{align}
 \| u \|_{X^{\frac{1}{4},\frac12+\eps}_{p}(I)}
\leq C_0 \|u(t_0) \|_{M^{2, p}_\frac{1}{4}}
\label{X7}
\end{align}

\noi
for some absolute constant $C_0 > 0$.
Then, 
from the Duhamel formula, 
Lemma \ref{LEM:lin} (with $ b = \frac 12+\eps$ and $b' = - \frac 12 + 2\eps$), 
\eqref{X6}, and \eqref{X7}, we obtain
\begin{align}
\begin{split}
\| u \|_{X^{s,\frac12+\eps}_{p}(I)}
& \les \| u(t_0) \|_{M^{2, p}_s} 
+   \dl^\eps
\| u \|_{X^{\frac{1}{4},\frac12+\eps}_{p}(I)}^2
\| u \|_{X^{s,\frac12+\eps}_{p}(I)}  \\
& \les \| u(t_0) \|_{M^{2, p}_s} 
+   \dl^\eps
\|u(t_0) \|_{M^{2, p}_\frac{1}{4}}^2 
\| u \|_{X^{s,\frac12+\eps}_{p}(I)}. 
\end{split}
\label{X8}
\end{align}

\noi
In particular, from \eqref{embed1} and \eqref{X8}, we conclude that 
there exists an absolute constant $C_1 > 0$ such that 
\begin{align}
\sup_{t \in [t_0, t_0 + \dl]}
\|u(t) \|_{M^{2, p}_s} 
\leq C_1 
\|u(t_0) \|_{M^{2, p}_s}
\label{X9}
\end{align}

\noi
for any $t_0 \in \R$.
Then, by iterating the local argument with
\eqref{X6a}, 
we conclude from \eqref{X9} that 
\begin{align*}
\sup_{t \in [0, T]} \| u(t) \|_{M^{2, p}_s} \leq C^{ (1 + \|u_0\|_{M^{2, p}_\frac{1}{4}})^\ta T} \|u_0\|_{M^{2, p}_s}
\end{align*}

\noi
for any $T > 0$.
This proves global well-posedness of \eqref{mKdV} in $M^{2, p}_s(\R)$
for $s \geq 1 - \frac 1p$.

\section{On the failure of local uniform continuity below $H^\frac{1}{4}(\R)$}
\label{SEC:ill}

In this section, we present the proof of Proposition \ref{PROP:ill}.
In particular, by  adapting the argument in \cite{KPV01} to the modulation space
setting, we prove the following statement.

\begin{lemma}\label{LEM:ill}
Suppose that $(s, p)$ satisfies one of the following conditions:
\textup{(i)}  $2 \leq p  \leq  \infty$ and $0 \leq  s < \frac 14$
or \textup{(ii)} $2 \leq p  <  \infty$ and 
 $ - \frac 1p <  s < 0$.
There exist two sequences $\{u_{0, n}\}_{n \in \N}$ and $\{\wt u_{0, n}\}_{n \in \N}$
in $\S(\R)$ 
such that 
\begin{itemize}

\item[\textup{(a)}]  
$u_{0, n}$ and $\wt u_{0, n}$ are uniformly bounded in $M^{2, p}_s(\R)$,

\smallskip

\item[\textup{(b)}]  
$\displaystyle \lim_{n \to \infty} 
\|u_{0, n}- \wt u_{0, n}\|_{M^{2, p}_s} = 0$,

\smallskip

\item[\textup{(c)}]  
Let $u_n$ and $\wt u_n$ be the solutions
to the focusing mKdV \eqref{mKdV} \textup{(}with the $+$ sign\textup{)} with initial data
$u_n|_{t = 0} =u_{0, n}$
and  
$\wt u_n|_{t = 0} =\wt u_{0, n}$, respectively.
Then, 
there exists $c > 0$ such that 
\begin{align*}
 \liminf_{n \to \infty}
\| u_n(T) - \wt u_n(T)\|_{ M^{2, p}_s} \geq c > 0
\end{align*}

\noi
for any $T > 0$.

\end{itemize}

\end{lemma}

In \cite{KPV01}, Kenig-Ponce-Vega
proved Lemma \ref{LEM:ill} 
for $p = 2$ 
by using 
explicit soliton solutions with parameters (see \eqref{un} below).
In the following, we use exactly the same explicit
soliton solutions to show an analogous instability
in the modulation space setting.

Let 
\begin{align}
\label{soli1}
Q(x) =   \sech (x).
\end{align}

\noi
Then, $f$ solves the ODE: $- Q + Q'' + 2Q^3 = 0$ and hence
\[
 - Q' + Q'''  + 6 Q^2 Q' =0.
\]

\noi
With  $Q_\ld (x) = \ld Q (\ld x)$, define $u_{N, \ld}$  by 
\begin{align}
\label{un}
u_{N,\ld} (x,t) = \frac{1}{\sqrt{6}} e^{it (  N^3 - 3N \ld^2) + iNx} Q_\ld (x+3N^2 t -\ld^2 t)
\end{align}

\noi
for  $N, \ld > 0$.
Then, it is easy to check that $u_{N,\ld}$  is a solution to \eqref{mKdV} 
with 
$u_{N,\ld}|_{t = 0}  = \frac{1}{\sqrt 6}e^{iN x} Q_\ld$
for any $N, \ld > 0$.
Recalling that 
\begin{align*}
\ft Q_\ld (\xi) = \ft Q \Big( \tfrac{\xi}{\ld} \Big) = \pi \sech \bigg( \frac{\pi \xi}{2\ld}\bigg), 
\end{align*}

\noi
we have
\begin{align}
\ft Q_\ld (\xi) \sim e^{- \frac{\pi |\xi|}{2\ld}}.
\label{Y2}
\end{align}

\noi
In particular, when $\ld \gg 1$, 
it follows from \eqref{un} that $\ft u_{N, \ld}(\xi, t)$ is highly concentrated around $|\xi|\sim N$.
See \eqref{Y4} below.

In the following, we first present the argument for $0 \leq s < \frac 14$.
We then discuss the case for $- \frac 1p < s < 0$ in Subsection \ref{SUBSEC:neg}.

\subsection{On the size of the soliton solutions}

Fix $2\leq p \leq \infty$ and  $0 \leq  s < \frac 14$.
Given $N \geq 1$, 
we consider two solutions $u_{N_1, \ld}$ and $u_{N_2, \ld}$
of the form \eqref{un}, 
where
\begin{align}
\ld = N^{-2s} \qquad \text{and}
\qquad N_1,N_2 =  N + O(1)
\label{para}
\end{align}

\noi
As we see below, we also impose
that  $|N_1 - N_2| \ll 1$.
Furthermore, fix 
  $\ta = \ta(s) > 0$  such that 
\begin{align}
4s - 1 + 2\ta < 0. 
\label{Y0}
\end{align}

In the following, we  estimate the $M^{2, p}_s$-norms
of $u_{N_j, \ld}$, $j = 1, 2$.
Noting that 
$|\ft u_{N_j, \ld}(\xi, t) | = |\ft u_{N_j, \ld}(\xi, 0) |$, 
the following computation holds uniformly in $t\in \R$.
We separately consider the contributions
from (i) $|\xi - N | \ll N^\ta$
and  
(ii) $|\xi - N | \ges N^\ta$.
Set 
\begin{align}
u^{(1)}_{N_j, \ld}
= \F_x^{-1} 
\big(\ind_{|\xi - N| \ll N^\ta} \cdot \ft u_{N_j, \ld}\big)
\qquad \text{and}\qquad
u^{(2)}_{N_j, \ld}= u_{N_j, \ld} - u^{(1)}_{N_j, \ld}.
\label{Y3}
\end{align}

We first consider (ii).
Note that when 
 $|\xi - N | \ges N^\ta$
 and 
$|\xi |\ges N$, we have
 $|\xi - N | \ges |\xi|^\ta$
 for small $\ta >0$.
Then, 
by separately considering the contribution
from $|\xi |\ll N$ and $|\xi |\ges N$, 
it follows from
\eqref{un}, 
\eqref{Y2},   \eqref{para}, and \eqref{Y0} that 
\begin{align}
\| u^{(2)}_{N_j, \ld} (t) \|_{M^{2,p}_s} 
& \sim  
 \bigg( \sum_{|n - N| \ges  N^\ta }
\jb{n}^{sp}
 e^{-\frac{p\pi }{2}  N^{2s} |n - N_j|} \bigg)^{\frac1p} \notag \\
& \les  
 \bigg( \sum_{|n| \ll N}
\jb{n}^{sp}
 e^{-\frac{p\pi }{2}  N^{\ta+ 2s} } \bigg)^{\frac1p} 
 + 
 \bigg( \sum_{|n| \ges N}
\jb{n}^{sp}
 e^{-\frac{p\pi }{2}  N^{2s} |n|^\ta} \bigg)^{\frac1p} \notag \\
& \les  e^{-c  N^{\ta+ 2s} }
\label{Y4}
\end{align}

\noi
since $\ta + 2s > 0$.
On the other hand, 
by a change of variables with \eqref{para}
and \eqref{Y2}, we have
\begin{align}
\| u^{(1)}_{N_j, \ld} (t) \|_{M^{2,p}_s} 
& \leq \| u^{(1)}_{N_j, \ld} (t) \|_{H^s} \notag \\
& \les   N^{s} \bigg(  \int_{|\xi - N|\ll N^\ta} | \ft Q_\ld (\xi -N_j)  |^2  d\xi \bigg)^{\frac{1}2} \notag \\
&  =    \bigg(  \int_{|\xi|\ll N^{\ta+2s}} e^{-\pi |\xi|} d\xi \bigg)^{\frac{1}2} \notag \\
& \sim 1.
\label{Y5}
\end{align}

\noi
By considering the contribution from $|\xi - N |\les 1$, 
we also see that 
\begin{align}
\| u^{(1)}_{N_j, \ld} (t) \|_{M^{2,p}_s} \ges 1.
\label{Y6}
\end{align}

\noi
Hence, from \eqref{Y3}, \eqref{Y4}, \eqref{Y5}, and \eqref{Y6}, 
we conclude that 
\begin{align}
\| u_{N_j, \ld} (t) \|_{M^{2,p}_s} \sim 1
\label{Y7}
\end{align}

\noi
for any $t \in \R$, independent of $N, N_1, N_2\geq 1$.

\subsection{On the difference of the soliton solutions}
\label{SUBSEC:diff}

When $t = 0$, we have the following upper bound
from \cite[(3.5)]{KPV01}:
\begin{align}
\begin{split}
\| u_{N_1,\ld} (0) - u_{N_2,\ld} (0) \|_{M^{2,p}_s}
& \leq 
\| u_{N_1,\ld} (0) - u_{N_2,\ld} (0) \|_{H^s}\\
& \les N^{2s}|N_1 - N_2|.
\end{split}
\label{diff}
\end{align}

Fix $T > 0$.
We establish a lower bound on the $M^{2, p}_s$-norm of the difference of $u_{N_j, \ld}(T)$.
In view of \eqref{Y3}
and \eqref{Y4}, it suffices to consider
$u^{(1)}_{N_1, \ld}(T) - u^{(1)}_{N_2, \ld}(T)$.
As in \cite{KPV01}, 
the main ingredient is separation of the physical supports 
of the soliton solutions $u_{N_j, \ld}$, $j = 1, 2$.
From~\eqref{un} with \eqref{soli1}, 
we see that $u_{N_j,\ld}(T)$  is concentrated on an interval 
of length $\sim \ld^{-1}$ centered at  $3N_j^2T - \ld^2 T$.
Note that these essential supports of
$u_{N_j,\ld}(T)$, $j = 1, 2$
 are disjoint,  provided that
\begin{align*}
N |N_1 -N_2| T \gg \ld^{-1} = N^{2s}.
\end{align*}

In our modulation space setting, however, 
we need to establish 
 separation of the physical supports 
of the frequency localized soliton solutions $\Pi_n u_{N_j, \ld}$, $j = 1, 2$.
From \eqref{Pi1}, 
there exists $\eta \in \S(\R)$ such that 
\[ |\Pi_n^2 u(x)| \leq  (|\eta| * |u|)(x) \]

\noi
for any $x \in \R$ and $n \in \Z$.
Then, from 
\eqref{un} and \eqref{soli1}, we have 
\begin{align}
\big|\jb{\Pi_n  & u_{N_1,\ld}(T), \Pi_n u_{N_2,\ld}(T)}_{L^2_x} \big|
 = \big|\jb{\Pi_n^2 u_{N_1,\ld}(t),  u_{N_2,\ld}(T)}_{L^2_x} \big|\notag \\
& \les \int_\R 
\bigg( \int_\R  |\eta(y)| Q_{\ld} (x-y + 3N^2_1 T -T\ld^2 ) dy \bigg) 
Q_{\ld} (x+ 3N^2_2 T-T\ld^2 )  dx\notag \\
& =  \int_\R 
\bigg( \int_\R  |\eta(y)| Q_\ld (x-y) dy \bigg)Q_{\ld} (x + 3(N^2_2-N^2_1) T)  dx\notag \\ 
& =   \iint_\R  |\eta(\ld^{-1} y)|  Q (x - y) Q (x + 3\ld(N^2_2-N^2_1) T) dy  dx\notag \\ 
& \les  \iint_{\R} \frac{1}{\jb{\ld^{-1} y }^K} e^{-|x-y|} e^{-|x + 3\ld(N^2_2-N^2_1) T|} dy  dx\notag \\ 
& \les 
\frac{1}{ N |N_1-N_2| T}
\label{Y8}
\end{align}

\noi
uniformly in  $n \in \Z$.

Given $N \gg 1$, choose $N_1, N_2 \sim N$ such that 
\begin{align}
|N_1 - N_2| \sim \frac{N^{2s - 1 +2 \ta}}{T}, 
\label{Y8a}
\end{align}

\noi
where $\ta > 0$ is as in \eqref{Y0}.
Thus, from  the triangle inequality,  \eqref{Y8}, 
Minkowski's inequality, and 
\eqref{Y6}
we have 
\begin{align}
\|u^{(1)}_{N_1, \ld}(T) & - u^{(1)}_{N_2, \ld}(T)\|_{M^{2, p}_s}^2
 \sim N^{2s} 
\bigg( \sum_{|n - N|\ll N^\ta  } \| \Pi_n u_{N_1, \ld}(T) 
- \Pi_n u_{N_2, \ld}(T) \|_{L^2_x}^p\bigg)^\frac{2}{p} \notag \\
& =  N^{2s} 
\bigg( \sum_{|n - N|\ll N^\ta  } \Big(\| \Pi_n u_{N_1, \ld}(T)\|_{L^2_x}^2  
+ \|\Pi_n u_{N_2, \ld}(T) \|_{L^2_x}^2  \notag\\
& \hphantom{XX}
- 2 \Re \jb{\Pi_n   u_{N_1,\ld}(t), \Pi_n u_{N_2,\ld}(t)}_{L^2_x} \Big)^\frac{p}{2}
 \bigg)^\frac{2}{p}\notag \\
& \ges
\| u^{(1)}_{N_1, \ld}(T)\|_{M^{2, p}_s}^2  
- N^{\frac{2}{p}\ta +2s} N^{-2\ta - 2s}\notag\\
& \ges 1 - N^{-2\ta (1 - \frac{1}{p})}
\sim 1
\label{Y9}
\end{align}

\noi
for any sufficiently large $N \gg 1$.
Hence, 
from \eqref{Y4}, \eqref{Y7},  and \eqref{Y9}, 
we conclude that 
\begin{align}
\|u_{N_1, \ld}(T)  - u_{N_2, \ld}(T)\|_{M^{2, p}_s}
\sim 1.
\label{Y10}
\end{align}

\noi
On the other hand, 
from \eqref{diff} and \eqref{Y8a} with \eqref{Y0}, 
we have
\begin{align}
\|u_{N_1, \ld}(0)  - u_{N_2, \ld}(0)\|_{M^{2, p}_s}
& \sim T^{-1} N^{4s-1+2 \ta} \notag\\
& \longrightarrow 0
\label{Y11}
\end{align}

\noi
by taking $N \to \infty$.
Finally, given $n \in \N$, 
let $N = 2^n$
and set 
$u_n = u_{N_1(n), \ld(n)}$
and $\wt u_n = u_{N_2(n), \ld(n)}$, 
Lemma \ref{LEM:ill} and hence Proposition \ref{PROP:ill}
follow from 
\eqref{Y7}, \eqref{Y10}, and~\eqref{Y11},
provided $2 \leq p \leq \infty$ and $0 \leq  s < \frac 14$.

\subsection{Failure of 
 uniform continuity in negative regularities}
\label{SUBSEC:neg}

In this subsection, we briefly consider the case $s < 0$.
With $\ld = N^{-2s}$ as in \eqref{para}, 
the estimate \eqref{Y7} is no longer true
and hence \eqref{Y10} fails in this case.

Fix $2 \leq p < \infty$ and $-\frac 1p <  s< 0$.
In the following, we use a new choice for the parameter~$\ld$:
\begin{align}
\ld = N^{-ps}
\label{Z1}
\end{align}

\noi
and let  $u_{N_j, \ld}$, $j = 1, 2$, be 
the solutions of the form \eqref{un} with this choice of $\ld$
(and $N_1, N_2 \sim N$).
We also choose new $\ta = \ta(s, p)>0$
such that 
\begin{align}
 -  ps < \ta < 1. 
\label{Z2}
\end{align}

\noi
This imposes the lower bound:  $s > -\frac 1p$.
Note that 
$|n - N|\ll N^\ta$ implies $|n| \sim N$
since $\ta < 1$.

By repeating the computation in \eqref{Y4}, we have
\begin{align}
\| u^{(2)}_{N_j, \ld} (t) \|_{M^{2,p}_s} 
& \les  e^{-c  N^{\ta+ ps} }
\label{Z3}
\end{align}

\noi
thanks to \eqref{Z2}.
On the other hand, 
from \eqref{un} and $Q_\ld(x) = \ld Q(\ld x)$, we have
\begin{align}
\| u^{(1)}_{N_j, \ld} (t) \|_{M^{2,p}_s}
& \sim N^s \Bigg( \sum_{|n - N|\ll N^\ta}
\bigg(\int_n^{n+1} |\ft Q_\ld(\xi - N_j)|^2 d\xi\bigg)^\frac{p}{2}\Bigg)^\frac{1}{p} \notag \\
& = N^s \ld^\frac{1}{2} 
\Bigg( \sum_{|n |\ll N^{\ta}}
\bigg(\int_\frac{n}{\ld}^{\frac{n+1}{\ld}} |\ft Q(\xi)|^2 d\xi\bigg)^\frac{p}{2}\Bigg)^\frac{1}{p} \notag \\
& \sim  N^s \ld^\frac{1}{p} 
\bigg( \sum_{|n |\ll N^{\ta} }
\big |\ft Q\big(\tfrac{n}{\ld}\big)\big|^p \ld^{-1}\bigg)^\frac{1}{p} \notag \\
\intertext{By the Riemann sum approximation with
\eqref{Z1} and \eqref{Z2},}
& \sim \| Q \|_{\FL^{0, p}} \sim 1, 
\label{Z5}
\end{align}

\noi
uniformly in large $N \gg 1$.
Hence, from \eqref{Y3}, \eqref{Z3}, and \eqref{Z5}, 
we conclude that 
\begin{align}
\| u_{N_j, \ld} (t) \|_{M^{2,p}_s}
\sim \| u_{N_j, \ld}^{(1)} (t) \|_{M^{2,p}_s} \sim 1
\label{Z6}
\end{align}

\noi
for any $t \in \R$, independent of $N, N_1, N_2\geq 1$.

Next, we estimate the difference of the soliton solutions as in Subsection \ref{SUBSEC:diff}.
A direct computation as in  \cite[(2.10)]{KPV01}
shows that 
\begin{align}
\begin{split}
\| u_{N_1,\ld} (0) - u_{N_2,\ld} (0) \|_{M^{2,p}_s}
& \les N^{s}\ld^{-\frac{1}{2}} |N_1 - N_2|.
\end{split}
\label{Z7}
\end{align}

\noi
In estimating the difference at time $T > 0$, 
we once again use the almost orthogonality
of the two soliton solutions, 
provided that 
\begin{align*}
N |N_1 -N_2| T \gg \ld^{-1} = N^{ps}.
\end{align*}

Given $N \gg 1$, choose $N_1, N_2 \sim N$ such that 
\begin{align}
|N_1 - N_2| \sim  \frac{N^{ps - 1 + \frac 32\ta }}{T}.
\label{Z9}
\end{align}

\noi
Then, by proceeding as in \eqref{Y9} with \eqref{Z6} and \eqref{Y8}
and choosing $\theta > - ps$ sufficiently close to $-ps$, we obtain
\begin{align}
\|u^{(1)}_{N_1, \ld}(T)  - u^{(1)}_{N_2, \ld}(T)\|_{M^{2, p}_s}^2
& \ges
\| u^{(1)}_{N_1, \ld}(T)\|_{M^{2, p}_s}^2  
- N^{\frac{2}{p}\ta +2s} N^{-ps- \frac 32\ta }\notag\\
& \ges 1
\label{Z10}
\end{align}

\noi
for all sufficiently large  $N \gg 1$.
Hence, 
from \eqref{Z3}, \eqref{Z6},  and \eqref{Z10}, 
we conclude that 
\begin{align*}
\|u_{N_1, \ld}(T)  - u_{N_2, \ld}(T)\|_{M^{2, p}_s}^2
\sim 1.
\end{align*}

\noi
On the other hand, 
from \eqref{Z7} with \eqref{Z1} and \eqref{Z9}, 
we have
\begin{align*}
\|u_{N_1, \ld}(0)  - u_{N_2, \ld}(0)\|_{M^{2, p}_s}^2
& \sim T^{-1} N^{s + \frac 32(\ta + ps) -1} \notag\\
& \longrightarrow 0
\end{align*}

\noi
by taking $N \to \infty$
since we chose  $\theta > - ps$ sufficiently close to $-ps$.
This completes the proof of 
Lemma \ref{LEM:ill} and hence Proposition \ref{PROP:ill}
when $2 \leq p < \infty$ and $-\frac 1p < s < 0$.

\begin{remark}\rm
Note that our parameter
 choices \eqref{para} for $s \geq 0$ and \eqref{Z1} for $s < 0$
 agree with those in \cite{KPV01} and \cite{Grun}, respectively.
In the following, we provide an intuitive explanation of our choices.
Given $f \in \S(\R)$, let $f_{N, \ld}(x) = \ld e^{iNx} f(\ld x)$.
When $s> 0$, $\ld = N^{-2s}$ in \eqref{para} tends to 0 as $N \to \infty$.
This implies that $\ft f_{N, \ld}$ is highly localized around $|\xi - N| \les \ld$. 
Namely, $\ft f_{N, \ld}$ is essentially supported in one interval $[N - \frac12, N + \frac 12)$,
in which case the $M^{2, p}_s$-norm of $\ft f_{N, \ld}$ 
reduces to its $H^s$-norm (which in turn can be reduced to the $L^2$-norm of $f$).
Therefore, the choice $\ld = N^{-2s}$ from the $H^s$-theory in \cite{KPV01} is appropriate
in this case.

On the other hand, when $s < 0$, 
$\ld = N^{-ps}$ tends to $\infty$ as $N \to \infty$.
Namely, the essential support of 
$\ft f_{N, \ld}$ spreads out as $N \to \infty$.
Then, arguing as in \eqref{Z5}, we see that 
the $M^{2, p}_s$-norm of $\ft f_{N, \ld}$
essentially reduces to the $\FL^{0, p}$-norm of $f$,
which shows that the choice $\ld = N^{-ps}$ from the $\FL^{s, p}$-theory in \cite{Grun} is appropriate
in this case.

\begin{ackno}\rm
T.O.~was supported by the ERC starting grant 
no.~637995 ``ProbDynDispEq''.

\end{ackno}

\end{remark}

\end{document}